\documentclass[11pt]{article}
\usepackage{amsmath, amssymb, amsthm}
\usepackage{graphicx}
\usepackage{color}
\usepackage{subfigure}
\usepackage{lineno}

\date{}
\begin{document}

{\LARGE \bf  
\begin{center}
Structural sensitivity of chaotic dynamics in\\ Hastings-Powell's model
\end{center}
}

\vspace*{1cm}

\centerline{\bf Indrajyoti Gaine$^a$, Swadesh Pal$^{b}$, Poulami Chatterjee$^{c}$, Malay Banerjee$^{a,}$\footnote{Corresponding author: malayb@iitk.ac.in}}

\vspace{0.5cm}

\centerline{ $^a$ Indian Institute of Technology Kanpur, Kanpur - 208016, India}

\centerline{ $^b$ MS2Discovery Interdisciplinary Research Institute, Wilfrid Laurier University,}

\centerline{75 University Ave W, Waterloo, N2L3C5, Ontario, Canada}

\centerline{ $^c$ Department of Mathematics, Jadavpur University, Kolkata, India} 

\vspace{1cm}

\begin{center}
{\bf Abstract}
\end{center}

The classical Hastings-Powell model is well known to exhibit chaotic dynamics in a three-species food chain. Chaotic dynamics appear through period-doubling bifurcation of stable coexistence limit cycle around an unstable interior equilibrium point. A specific choice of parameter value leads to a situation where the chaotic attractor disappears through a collision with an unstable limit cycle. As a result, the top predator goes to extinction. Here we explore the structural sensitivity of this phenomenon by replacing the Holling type II functional responses with Ivlev functional responses. Here we prove the existence of two Hopf-bifurcation thresholds and numerically detect the existence of an unstable limit cycle. The model with Ivlev functional responses does not indicate any possibility of extinction of the top predator. Further, the choice of functional responses depicts a significantly different picture of the coexistence of the three species involved with the model.

\vspace{1.0cm}



\baselineskip 0.25in

\begin{center}
{\LARGE\bf }
\end{center}

\section{Introduction}

 In an ecosystem, various species interact with one another in a variety of ways, such as through mutualism, competition, and predation. Mathematical modelling helps to capture such phenomena and predicts their dynamics for the long term based on the current state of the population and knowledge of relevant ecological processes. Researchers have been using different types of mathematical modelling approaches to take care of these phenomena, e.g., ordinary differential equations (ODEs), partial differential equations, etc. Different mathematical forms are also used in the model to incorporate the underlying interactions. ODE models provide a framework to describe the dynamics of population growth over time, and these are commonly studied for single-species, two-species, and multi-species interactions within ecological systems. The equilibrium points and their stabilities are the two important factors investigated in the ODE model, which helps in understanding the local behaviour of an ecological system, in particular, the co-existing equilibrium point(s) since it provides information about the coexistence of species.

 Different types of coexistence behaviours occur in an ecological model, and they change through local and global bifurcations. For instance, the stable coexistence of a system shows constant dynamics over a long time; however, a regular oscillation gives synchronised periodic dynamics. These two coexistences can switch in between through a supercritical Hopf bifurcation \cite{perko2013}. Ecological systems are inherently nonlinear and can exhibit unexpected fluctuations and irregular oscillations, called chaos. It depends on the background parameters and initial conditions, small changes which lead to exponential deviations. It is known that the single- and two-species autonomous ODE models do not generate irregular oscillatory coexistence; however, models with three or more species can exhibit a wide range of oscillatory solutions, including quasi-periodic and chaotic. These two types of oscillations often explain the co-existence of all the species with varying amplitude. There are many kinds of complex dynamics in a population model, such as multiple attractors \cite{tang2003}, catastrophic transitions \cite{hilker2010}, sub-harmonics of various periods \cite{schaffer2001}, cascades of period doubling \cite{previte2013}, and strange attractors \cite{klebanoff1994A}.

 According to Bazykin \cite{bazykin1998}, a two-species model has stable equilibrium and stable cycles separated by unstable cycles. Though this model produces complicated dynamics, it does not produce complex oscillations like chaos, as shown in experimental research. In carefully monitored laboratory trials, cultures of our beetles (Tribolium castaneum) experience bifurcations in their dynamics when the demographic factors change, including a specific route to chaos \cite{cushing2001}. Furthermore, a long-term experiment on a complex food web shows chaos in the ecology \cite{tanabe2005}. The chaotic dynamics in food-chain models emphasize the ecological systems' innate complexity and how challenging it is to forecast their behaviour. Researchers have observed chaotic dynamics in continuous systems, such as the two-prey-one-predator model \cite{gilpin1979}, one-prey-two-predator model \cite{klebanoff1994B}, two-sex model \cite{caswell1986}, three-species food chain models \cite{hastings1991, mccann1994A, mccann1994B}.

 Chaos appears in mathematical models through periodic doubling \cite{may1976A,stone1993}. On the other hand, the disappearance of chaos does not always follow the same pattern \cite{boer2001}. In general, when chaos is suppressed, it is either through global bifurcations or crises. In the case of global bifurcations, the basin boundary collision happens between several basins of attraction (areas of the phase space where starting circumstances converge to a certain attractor or behaviour). The crisis happens between the collision of a chaotic attractor and a co-existing unstable fixed point or periodic orbit \cite{boer2001}. As a result, the system shifts to a sudden qualitative change in dynamics and exhibits constant or periodic behaviour for both cases following the disappearance.

 The systems' dynamics depend on the functional form or the parameter values involved in it. For instance, two types of population growth are considered in the literature: exponential and logistic growth. The growth rate for exponential growth is density-independent, and the logistic growth is density-dependent, the most common growth rate considered in ecological models. In addition, the functional response is also responsible for the resultant dynamics as it links between two connected trophic levels and plays a crucial role in shaping the dynamics of the system. Primarily, the functional response can be classified into two groups: prey-dependent \cite{holling1959}, and prey-predator-dependent \cite{beddington1975, deangelis1975}. Sometimes, two different types of functional response can show different dynamics, and it is based upon predatory interactions \cite{oaten1975, cantrell2001}.

Researchers have been exploring how different parametrization has significantly altered the model dynamics. For the past decades, they have focused on how different functional forms with similar geometrical shapes and likewise properties influence the model dynamics \cite{myerscough1996, gross2004, wood1999, fussmann2005}. A small change in functional form, here the functional response, can cause significant consequences in the dynamics of the model, referred to as ‘structural sensitivity' \cite{adamson2014}. It has been studied on different ecological models, e.g., zooplankton feeding on multiple resources \cite{gentleman2003}, nutrient–phytoplankton–zooplankton model \cite{gentleman2008}, and complex marine ecosystem \cite{anderson2010}, etc. Different approaches also have been applied to study structural sensitivity, such as a probabilistic viewpoint \cite{aldebert2018}, the discrete-time system \cite{rana2020}, a statistical viewpoint \cite{aldebert2016A}, and many more. It is shown that if functional responses of a sub-model fit in the area of the co-existing stable state equilibrium point, then the structural sensitivity reduces \cite{wyse2022}. Furthermore, researchers have shown that a small continuous change in the functional form leads to significant alterations in bifurcations from a three-dimensional point of view \cite{adamson2013,aldebert2019}.

 In ecology, different functional forms are used in mathematical models to incorporate the same type of physical phenomena. The primary focus of this work is to see the behavioural change in the dynamics by varying their forms, here the functional responses, along with the background parameters in terms of structural sensitivity. We have chosen an appropriate parameterization of these functional responses, which follows the assumptions: zero at zero, strictly monotonically increasing, finite horizontal asymptote, and concave down \cite{holling1959}. These assumptions help to understand the threshold explicitly behind a regime shift of systems' dynamics. They are needed because sometimes, a model can predict many results without any functional forms, but comparing or understanding their qualitative characteristics may become challenging. In \cite{seo2018}, authors have considered a two-species food chain model and studied local and global dynamics in the presence of different functional responses. Here, we extend this idea into a three-species model and study the sensitivity of chaotic dynamics. For this, we consider Hasting's model \cite{hastings1991}, which features a number of intriguing dynamics, including the "tea-cup" dynamics and the appearance of chaos through periodic orbit. In order to examine the model's structural stability, we take the parameter configurations in such a way that they look similar and consider the bifurcation parameter that is independent of these functional responses.

 The organization of the paper is as follows. The main structure of the model and the conditions for the functional response is given in Sect. \ref{MathM}. The model can have multiple equilibrium points depending on the parametric conditions. All the possibilities for the equilibrium points and their stability criteria are discussed in Sect. \ref{EqPtSt}. The model changes its dynamics by changing the bifurcation parameter (the mortality rate of the top predator), keeping others fixed, and this happens due to different types of bifurcations such as saddle-node, transcritical, and Hopf. The analytical conditions for these bifurcations in terms of the general functional form are presented in Sect. \ref{TempBif}. We have validated our theoretical findings by choosing the appropriate parametric values in Sect. \ref{NumRes}.

\section{Mathematical Model}{\label{MathM}}

We first consider a general form of a three-species model \cite{hastings1991} with logistic growth in the prey population and linear death rates ($d_{1}$ and $d_{2}$) for the intermediate and top predators as:
\begin{subequations}\label{genmodel}
\begin{align}
\frac{dx}{dt} &= x-x^2-f_{1}(x)y,\\
\frac{dy}{dt} &= f_{1}(x)y-d_1 y- f_{2}(y)z,\\
\frac{dz}{dt} &= f_{2}(y)z-d_2z,
\end{align}
\end{subequations}
where $f_{1}(\cdot)$ and $f_{2}(\cdot)$ are called the functional responses, and each $f_{i}(i=1,2)$ satisfies the conditions: (I) zero at zero, i.e., $f_{i}(0) = 0$; (II) monotone increasing, i.e., $f_{i}'(u) > 0~\forall u \geq 0$; (III) finite horizontal asymptote, i.e., $\lim_{u\rightarrow\infty}f_{i}(u) =  f_{i}^{\infty} < \infty$; (IV) concave down, i.e., $f_{i}''(u) < 0~\forall u \geq 0$. Researchers have been using different types of functional responses in ecological models to incorporate resource-consumer relations. Most of these functional responses satisfy only the first three conditions. But, we have added a condition (condition (IV)) for selecting those functional responses which have negative curvature. This extra condition gives us a more zoomed picture in the sensitivity analysis of the dynamics of the ecological model (\ref{genmodel}). In \cite{klebanoff1994A}, authors studied the model (\ref{genmodel}) by considering both the functional responses as Holling type II with different parametric setups. In this work, theoretically and numerically, we study the dynamics of the model for different mathematical forms of functional responses, particularly Holling type II and Ivlev, and compare their results qualitatively and quantitatively. We choose the expressions of the functions $f_{1}$ and $f_{2}$ for Holling type II functional responses as $f_{1}(x) = a_{1}x/(1+b_{1}x)$ and $f_{2}(y) = a_{2}y/(1+b_{2}y)$, and for Ivlev functional responses as $f_{1}(x) = \bar{a}_{1}(1-e^{-\bar{b}_1x})$ and $f_{2}(y) = \bar{a}_2(1-e^{-\bar{b}_2y})$.

\section{Equilibrium points and their stabilities}{\label{EqPtSt}}

Different approaches have been applied in the temporal model to find the nature of its solution(s). One of them is finding the equilibrium points and their stabilities. If an equilibrium point is locally stable, then the solution converges to itself, for the initial conditions lie in a certain neighbourhood around it. The linear stability analysis around the equilibrium point helps to find this type of local behaviour of the solution to the system. Before going to the stability, we first find the possible equilibrium points for the system (\ref{genmodel}), and these can be found by solving the system (\ref{genmodel}) with vanishing all derivatives, i.e., the solutions of the following algebraic equations:
\begin{subequations}\label{algeb_equn}
\begin{align}
x-x^2-f_{1}(x)y &= 0,\\
f_{1}(x)y-d_1 y- f_{2}(y)z &= 0,\\
f_{2}(y)z-d_2 z &= 0.\label{algeb_equnc}
\end{align}
\end{subequations}
The equilibrium points depend on the functional responses as they are involved on the left-hand side of the system (\ref{algeb_equn}). We use analytical and geometric approaches to find the equilibrium points.

The system (\ref{algeb_equn}) may have negative solutions, but we do not consider those as they correspond to negative densities, which are not feasible. From the expression of the algebraic system (\ref{algeb_equn}), we see that $E_{0} = (0,0,0)$ is the trivial solution, and it is an equilibrium point for the model (\ref{genmodel}). Furthermore, $E_{1} = (1,0,0)$ is also a solution to the system (\ref{algeb_equn}), which is an axial equilibrium point of the system (\ref{genmodel}). These two equilibrium points are independent of all possible forms of functional responses. There is a possibility of having a boundary equilibrium point of the system (\ref{genmodel}) of the form $E_b=(x_b,y_b,0)$ where $x_b=f_{1}^{-1}(d_1)$ and $y_b=x_b(1-x_b)/d_{1}$. According to our assumptions, the function $f_{1}$ is continuously strictly increasing and bounded above, so there may exist a unique solution of the equation $f_{1}(x) = d_{1}$ if $d_{1}<f_{1}^{\infty}$. In addition, the positivity condition for $y_{b}$ gives the feasibility of the equilibrium point $E_{b}$, which is $x_{b}<1$. This implies that $f_{1}^{-1}(d_1) <1$, which further implies $d_{1}<f_{1}(1)$. Combining all the conditions, the existence of the boundary equilibrium of the form $E_b=(x_b,y_b,0)$ where $x_{b}>0$ and $y_{b}>0$ requires the condition $d_{1} < \min\{f_{1}(1), f_{1}^{\infty}\}$. It is the only boundary equilibrium point that exists for the system (\ref{genmodel}). The considered system does not have any boundary equilibrium point on $xz$-plane where $x$ and $z$ are both non-zero. This is because $y=0$ holds on $xz$-plane, which further implies $z=0$ from (\ref{algeb_equnc}) and then we arrive at either at $E_{0}$ or $E_{1}$. Furthermore, the system does not have any boundary equilibrium on $yz$-plane where $y$ and $z$ are both non-zero as it corresponds to $x=0$, which gives $d_1y+d_2z = 0$, and it is not possible.

Now we discuss all possible numbers of interior equilibrium points for the system (\ref{genmodel}). The analytical and graphical approaches help us find such equilibrium points' existence and uniqueness. Before going into these, we first rearrange the algebraic equations (\ref{algeb_equn}) to satisfy an interior equilibrium point $E_{*} = (x_{*},y_{*},z_{*})$ as follows:
\begin{subequations}\label{re_algeb_equn}
\begin{align}
1-x &= \tilde{f}_{1}(x)y, \\
f_{1}(x) - d_1 &= \tilde{f}_{2}(y)z,\\
f_{2}(y) -d_2 &= 0,
\end{align}
\end{subequations}
where $\tilde{f}_{1}(x) = f_{1}(x)/x$ and $\tilde{f}_{2}(y) = f_{2}(y)/y$. These functions $\tilde{f}_{1}(x)$ and $\tilde{f}_{2}(y)$ are well-defined as $x$ and $y$ both are positive for the interior equilibrium point, and they are both positive, follow from the definitions of $f_{1}$ and $f_{2}$.

Since $f_{2}(y)$ is a strictly increasing function and has an upper limit $f_{2}^{\infty}$ (by the condition (II)), the third equation of (\ref{re_algeb_equn}) can have a unique positive root $y = y_{*}$ if $d_{2}<f_{2}^{\infty}$. Substituting this $y_{*}$ into the first equation of (\ref{re_algeb_equn}) gives an algebraic equation in terms of $x$, and its number of feasible solutions gives the maximum number of interior equilibrium points possible for the model. In this case, the solutions are the points of intersections between the line $y = 1-x$ and the curve $y = y_{*}\tilde{f}_{1}(x)$ in the first quadrant. Indeed, the points on the line $y = 1-x$ in the first quadrant satisfy $x<1$. The number of points of intersections depends on the characteristics of the function $\tilde{f}_{1}(x)$. Here, we assume some properties on $\tilde{f}_{i}$'s ($i = 1,2$) for finding the possible number of 
intersections: (a) $\lim_{u\rightarrow 0}\tilde{f}_{i}(u) > 0$, (b) $\tilde{f}_{i}'(u)< 0~\forall u \geq 0$, (c) $\tilde{f}_{i}''(u) > 0~\forall u \geq 0$, and (d) $\lim_{u\rightarrow\infty}\tilde{f}_{i}(u) = 0$.

\begin{figure}[h!]
\centering
\mbox{\subfigure[]{\includegraphics[width=5.5cm]{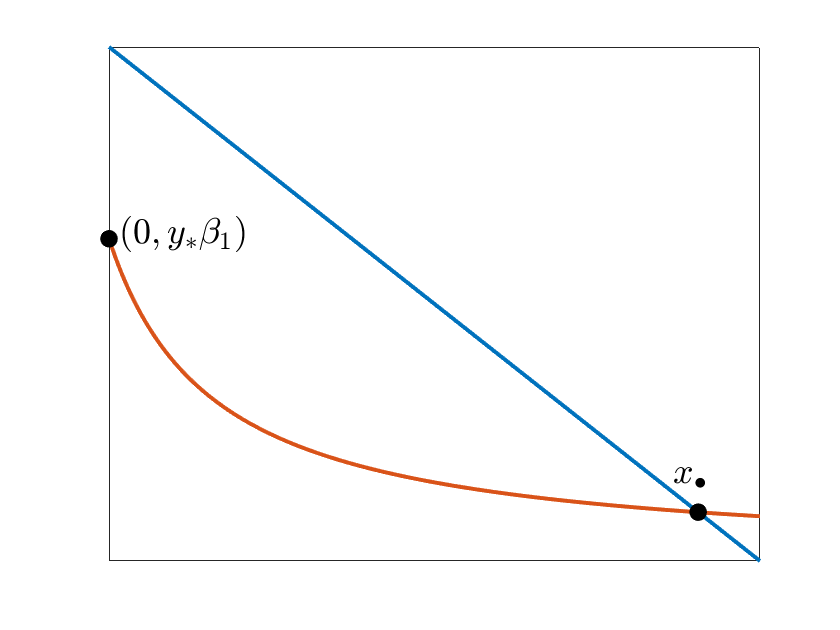}}
\subfigure[]{\includegraphics[width=5.5cm]{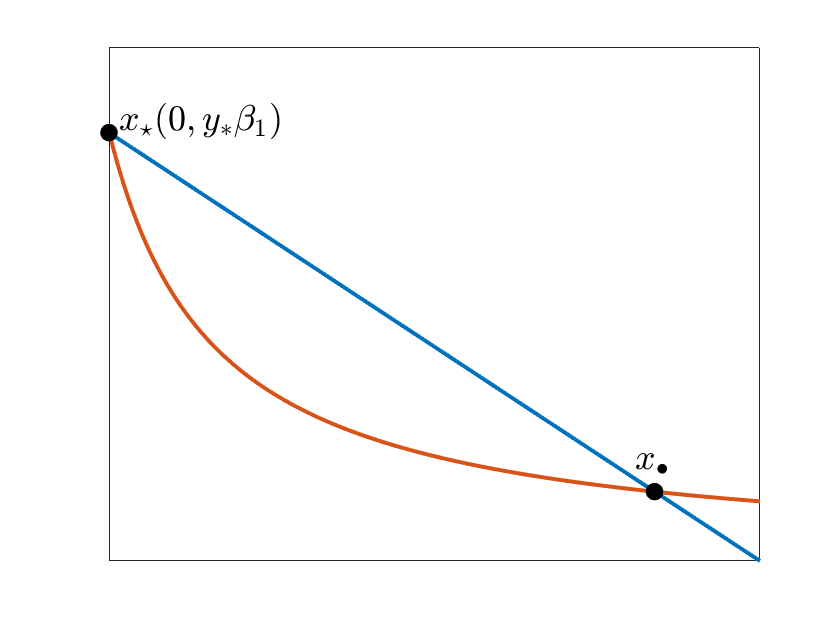}}}
\caption{Illustrations for the possible number of positive roots of the equation $1-x-y_*\tilde{f}_1(x)=0$ for $y_*\beta_1 \leq 1$ in $[0,1]$: (a) one root for $y_*\beta_1<1$ and (b) two roots for $y_*\beta_1=1$.}
\label{fig:figinterior2}
\end{figure}

\begin{figure}[h!]
\centering
\mbox{\subfigure[]{\includegraphics[width=5.5cm]{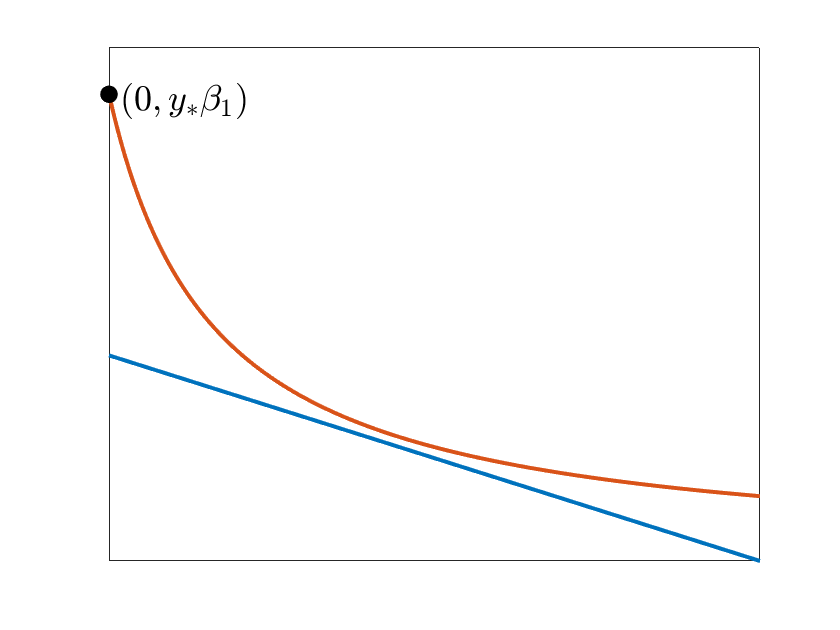}}
\subfigure[]{\includegraphics[width=5.5cm]{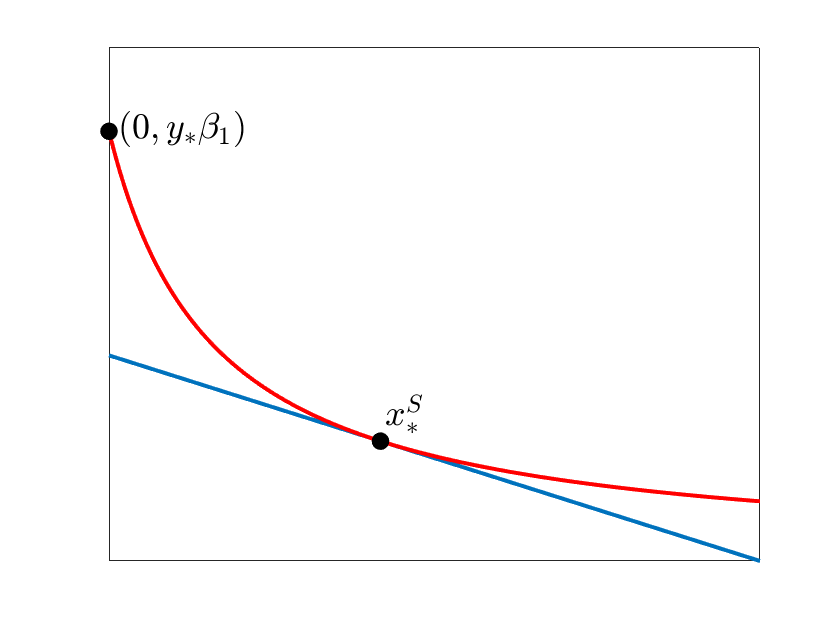}}
\subfigure[]{\includegraphics[width=5.5cm]{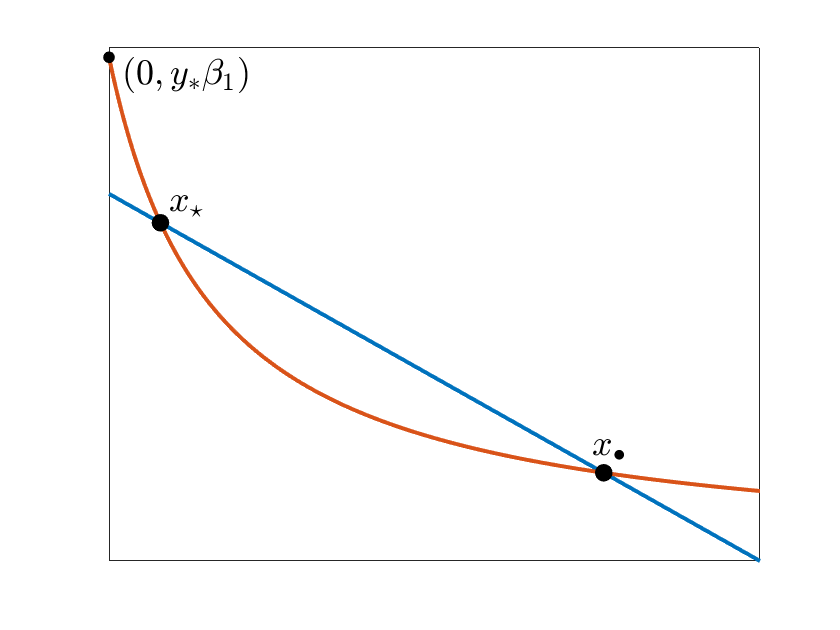}}}
\caption{Illustrations for the possible number of positive roots of the equation $1-x-y_*\tilde{f}_1(x)=0$ for $y_*\beta_1>1$ in $[0,1]$: (a) no root, (b) exactly one root, and (c) two roots.}
\label{fig:figinterior}
\end{figure}

 We set $\beta_{1} = \lim_{x\rightarrow 0}\tilde{f}_{1}(x)$ which is positive by the assumption (a). Also, the function $\tilde{f}_{1}$ is decreasing and convex by the assumptions (b) and (c). Therefore, if $y_{*}\beta_{1}<1$, then there exists a unique point of intersection between the curve $y = y_{*}\tilde{f}_{1}(x)$ and the line $y = 1-x$ [see Fig. \ref{fig:figinterior2}(a)]. On the other hand, for $y_{*}\beta_{1}\geq 1$, there may exists at most two intersections between the line $y = 1-x$ and the curve $y = y_{*}\tilde{f}_{1}(x)$ [see Figs. \ref{fig:figinterior2}(b) and \ref{fig:figinterior}]. In particular, the equation $1-x = y_{*}\tilde{f}_{1}(x)$ has two solutions for $x$ when $y_{*}\beta_{1} = 1$  [see Fig. \ref{fig:figinterior2}(b)]. Not only this, there could be two solutions to the equation $1-x = y_{*}\tilde{f}_{1}(x)$ for $y_{*}\beta_{1} > 1$. Furthermore, this equation can exist as a unique solution $x = x_{*}$, and in this case, the above-mentioned curve and line share the common tangent at $x = x_{*}$. When this occurs, the condition $y_{*}\tilde{f}_{1}'(x_{*}) = -1$ is satisfied and is referred to as a saddle-node bifurcation point at $d_{2} = d_{2}^{S}$, and is covered in greater detail later in this section.

For $d_{2}>d_{2}^{S}$, any solution $y_{*}$ for the equation $f_{2}(y) = d_2$ satisfies the inequality $y_{*}>y_{*}^{S}$ as $f_{2}(y)$ is an increasing function. Now, for this $y_{*}$, the equation $1-x = y_{*}\tilde{f}_{1}(x)$ does not have any solution because $1-x = y_{*}^{S}\tilde{f}_{1}(x)< y_{*}\tilde{f}_{1}(x)$ and $\tilde{f}_{1}(x)> 0~\forall x>0$. On the other hand, for $d_{2}<d_{2}^{S}$, the solution $y_{*}$ for the equation $f_{2}(y) = d_2$ satisfies the inequality $y_{*}<y_{*}^{S}$ and the equation $1-x = y_{*}\tilde{f}_{1}(x)$ exists two positive solutions for $x$ as $\tilde{f}_{1}(x)$ is convex and positive for $x>0$. In this case, we denote two solutions as $x_{\star}$ and $x_{\bullet}$, and in general, at these points $\tilde{f}_{1}(x)$ satisfies the conditions $y_{*}\tilde{f}_{1}'(x_{\star}) < -1$ and $y_{*}\tilde{f}_{1}'(x_{\bullet}) > -1$.

 Substituting the solution $x=x_{*}$ (whenever it exists) for the equation $1-x = y_{*}\tilde{f}_{1}(x)$ into the second equation of (\ref{re_algeb_equn}), one obtain $z_{*} = (f_{1}(x_{*})-d_{1})/\tilde{f}_{2}(y_{*})$. Depending on the number of solutions for $x$ for the equation $1-x = y_{*}\tilde{f}_{1}(x)$, one could get the same number of solutions for $z_{*}$. But, for each cases, the condition $f_{1}(x_{*})> d_{1}$ has to be satisfied for the feasibility of $z_{*}$. Suppose, two feasible solutions $z_{\star}$ and $z_{\bullet}$ exist for $x_{\star}$ and $x_{\bullet}$, respectively, and we denote the corresponding equilibrium points as $E_{\star} = (x_{\star},y_{\star},z_{\star})$ and $E_{\bullet} = (x_{\bullet},y_{\bullet},z_{\bullet})$ with $y_{*}=y_{\star} = y_{\bullet}$. Therefore, the inequality $y_{*}\tilde{f}_{1}'(x_{\star}) < -1$ implies $1-2x_{\star}-y_{\star}f_1'(x_{\star})>0$ and the other inequality $y_{*}\tilde{f}_{1}'(x_{\bullet}) > -1$ implies $1-2x_{\bullet}-y_{\bullet}f_{1}'(x_{\bullet})<0$. This further implies that both the interior equilibrium point satisfies $1-2x_{*}-y_{*}f_{1}'(x_{*})\neq0$ whenever $d_2<d_2^S$. It may happen that for some $x_{*}$, the condition $f_{1}(x_{*}) = d_{1}$ is satisfied, and in this case, the non-trivial equilibrium point coincides with the boundary equilibrium point $E_{b}$ on the $xy$-plane. We prove this bifurcation as a transcritical bifurcation, and it is discussed later in this section. 

 Now we investigate the model's local stability at various equilibrium points. This can be found by studying the eigenvalues of the Jacobian matrix at each equilibria. We start with the trivial equilibrium point $E_0=(0,0,0)$. The Jacobian matrix at $E_0=(0,0,0)$ is 
\begin{eqnarray*}
    \mathbf{J}_{E_{0}} = \left[
    \begin{array}{ccc}
    1 & 0 & 0 \\
    0 & -d_1 & 0 \\
    0 & 0 & -d_2 \\
    \end{array}\right].
\end{eqnarray*}
The eigenvalues of this matrix are $1$, $-d_1$, and $-d_2$. This shows that the Jacobian matrix has one positive and two negative eigenvalues. Therefore, the trivial equilibrium point $E_0$ is a saddle point, and it has a $2$-dimensional stable manifold and a $1$-dimensional unstable manifold, represented as $W^{s}(E_0)$ and $W^{u}(E_0)$, respectively. The stability of $E_0$ is independent of parametric restrictions. The Jacobian matrix evaluated at the axial equilibrium point $E_1=(1,0,0)$ is:
\begin{eqnarray*}
    \mathbf{J}_{E_{1}} = \left[
    \begin{array}{ccc}
    -1 & -f_{1}(1) & 0 \\
    0 & f_{1}(1)-d_1 & 0 \\
    0 & 0 & -d_2 \\
    \end{array}\right],
\end{eqnarray*}
which has the eigenvalues $-1$, $f_{1}(1)-d_{1}$, and $-d_{2}$. In this case, the Jacobian matrix has two negative eigenvalues for all parameter values, and the sign of the remaining one depends on the value of $f_{1}(1)$. If $f_{1}(1)>d_{1}$, then $E_1$ is a saddle point. On the other hand, for $f_{1}(1)<d_{1}$, $E_1$ is locally asymptotically stable with $3$-dimensional stable manifold. Similarly, we obtain the Jacobian matrix at the boundary equilibrium point $E_b=(x_b,y_b,0)$ as:
\begin{eqnarray*}
    \mathbf{J}_{E_{b}} = \left[
    \begin{array}{ccc}
    1-2x_b-y_bf_{1}'(x_b) & -f_{1}(x_b) & 0 \\
    y_bf_{1}'(x_b) & 0 & -f_{2}(y_b) \\
    0 & 0 & f_{2}(y_b)-d_2 \\
    \end{array}\right].
\end{eqnarray*}
This Jacobian has an eigenvalue $f_{2}(y_b)-d_2$, and the other two eigenvalues are the eigenvalues of a matrix having trace $1-2x_b-y_bf_{1}'(x_b)$ and determinant $y_bf_{1}(x_b)f_{1}'(x_b)$. Now, based on the properties of the functional responses $f_{i}$ ($i=1,2$), we can conclude that the boundary equilibrium point $E_b=(x_b,y_b,0)$ is locally asymptotically stable if $f_{2}(y_b)<d_2$ and $1-2x_b-y_bf_{1}'(x_b)<0$, and unstable for the other circumstance. For the unstable case, the dimension of its stable (unstable) manifold is determined by the number of negative (positive) eigenvalues of the Jacobian matrix $\mathbf{J}_{E_{b}}$.

 Finally, the Jacobian matrix evaluated at a typical interior equilibrium point $E_*=(x_*,y_*,z_*)$ is given by:
\begin{eqnarray*}
    \mathbf{J}_{E_{*}} = \left[
    \begin{array}{ccc}
    1-2x_*-y_*f_{1}'(x_*) & -f_{1}(x_*) & 0 \\
    y_*f_{1}'(x_*) & f_{1}(x_*)-d_1-z_*f_{2}'(y_*) & -d_{2} \\
    0 & z_*f_{2}'(y_*) & 0 \\
    \end{array}\right].
\end{eqnarray*}
The eigenvalues for this Jacobian matrix implicitly depend on the functional responses $f_{i}$ ($i=1,2$) and the equilibrium point $E_*=(x_*,y_*,z_*)$, and they are the solutions for the characteristics equation  $$\lambda^3+P_2\lambda^2+ P_1\lambda+P_0=0,$$
where $P_2=-(1-2x_*-y_*f_{1}'(x_*)+f_{1}(x_*)-d_1-z_*f_{2}'(y_*))$,
$P_1=(1-2x_*-y_*f_{1}'(x_*))(f_{1}(x_*)-d_1-z_*f_{2}'(y_*))+y_*f_{1}(x_*)f_{1}'(x_*)+d_{2}z_{*}f_{2}'(y_{*})$, and $P_0=-d_{2}z_{*}f_{2}'(y_{*})(1-2x_{*}-y_{*}f_{1}'(x_{*}))$. After applying the Routh-Hurwitz criteria, we can say the interior equilibrium point $E_*=(x_*,y_*,z_*)$ is stable if $P_2>0$, $P_0>0$ and $P_{1}P_{2} > P_{0}$ hold. 

\section{Temporal bifurcations}{\label{TempBif}}

 In this section, we discuss different temporal bifurcations for the system (\ref{genmodel}) where it changes the dynamics qualitatively and quantitatively through the bifurcation points, in particular, saddle-node, transcritical and Hopf bifurcations. As previously stated, we find all the bifurcation thresholds in terms of the bifurcation parameter $d_{2}$ along with the transversality conditions. In addition, we follow Sotomayor's theorem to verify the transversality conditions for these temporal bifurcations \cite{perko2013} with the following notation:
$$\mathbf{F}((x,y,z);d_2) \equiv \left[\begin{array}{c}
x-x^2-f_{1}(x)y \\
f_{1}(x)y-d_{1}y-f_{2}(y)z \\
f_{2}(y)z -d_{2}z\\
\end{array}\right].$$

\subsection{Saddle-node bifurcation}

 As we have mentioned in the previous section, the curve $y = y_{*}\tilde{f}_{1}(x)$ and the line $y = 1-x$ intersects at exactly one point in the first quadrant, and we have denoted the value for $d_{2}$ for which it occurs as $d_{2}^{S}$. Let us denoted that unique non-trivial equilibrium as $E_{*}^{S} = (x_{*}^{S},y_{*}^{S},z_{*}^{S})$. It is obvious that $d_{2}^{S}<f_{2}^{\infty}$, otherwise the solution $y= y_{*}^{S}$ does not exist for the third equation of (\ref{re_algeb_equn}). Furthermore, at $d_{2} = d_{2}^{S}$, the curve $y = y_{*}^{S}\tilde{f}_{1}(x)$ and the line $y = 1-x$ shares the common tangent, which gives the condition $y_{*}^{S}\tilde{f}_{1}'(x_{*}^{S}) = -1$. In this case, the Jacobian matrix at the non-trivial equilibrium point $E_{*}^{S}$ can be written as 
\begin{eqnarray*}
    \mathbf{J}_{E_{*}^{S}} = \left[
    \begin{array}{ccc}
    1-2x_{*}^{S}-y_{*}^{S}f_{1}'(x_{*}^{S}) & -f_{1}(x_{*}^{S}) & 0 \\
    y_{*}^{S}f_{1}'(x_{*}^{S}) & f_{1}(x_{*}^{S})-d_{1}-z_{*}^{S}f_{2}'(y_{*}^{S}) & -f_{2}(y_{*}^{S}) \\
    0 & z_{*}^{S}f_{2}'(y_{*}^{S}) & 0 \\
    \end{array}\right].
\end{eqnarray*}
Now, after simplify the condition $y_{*}^{S}\tilde{f}_{1}'(x_{*}^{S}) = -1$, one could obtain $1-2x_{*}^{S}-y_{*}^{S}f_{1}'(x_{*}^{S}) = 0$, which causes the determinant of the Jacobian matrix $\mathbf{J}_{E_{*}^{S}}$ to be zero. This implies that the Jacobian matrix $\mathbf{J}_{E_{*}^{S}}$ has a zero eigenvalue at $d_{2}^{S}$, which is simple, as the sum of the product of its eigenvalues is $y_{*}^{S}f_{1}'(x_{*}^{S})f_{1}(x_{*}^{S})+z_{*}^{S}f_{2}(y_{*}^{S})f_{2}'(y_{*}^{S})>0$. Therefore, $E_{*}^{S}$ is a non-hyperbolic equilibrium point. The eigenvectors corresponding to the zero eigenvalues of $\mathbf{J}_{E_{*}^{S}}$ and $[\mathbf{J}_{E_{*}^{S}}]^{T}$ are given by $$\mathbf{V} = \left[\begin{array}{c}
1 \\
0 \\
\frac{y_{*}^{S}f_{1}'(x_{*}^{S})}{f_{2}(y_{*}^{S})}\\
\end{array}\right] ~\mbox{and}~ \mathbf{W} = \left[\begin{array}{c}
1 \\
0 \\
\frac{f_{1}(x_{*}^{S})}{z_{*}^{S}f_{2}'(y_{*}^{S})}\\
\end{array}\right],$$ respectively. Furthermore, one finds the transversality conditions \cite{perko2013} as follows:
$$\mathbf{W}^T\mathbf{F}_{d_{2}}((x_{*}^{S},y_{*}^{S},z_{*}^{S});d_{2}^{S}) = \frac{-f_{1}(x_{*}^{S})}{f_{2}'(y_{*}^{S})} \neq  0$$
$$\mbox{and}~\mathbf{W}^T[D^2\mathbf{F}((x_{*}^{S},y_{*}^{S},z_{*}^{S});d_{2}^{S})(\mathbf{V},\mathbf{V})] = -2-f_{1}''(x_{*}^{S})y_{*}^{S} \neq 0.$$ 
This implies that the system undergoes a non-degenerate saddle-node bifurcation at $d_2 = d_{2}^{S}$.

\subsection{Transcritical bifurcation}

 Here we prove that the interior equilibrium point $E_{\star}$ can be generated from the boundary equilibrium point $E_{b}$ through a transcritical bifurcation. As discussed earlier, the interior equilibrium point coincides with the boundary equilibrium point when the solution $x_{*}$ of the equation $1-x-y_{*}\tilde{f}_{1}(x) = 0$ satisfies the condition $f_{1}(x_{*}) = d_{1}$. Assume this occurs at $d_2 = d_{2}^{T}$; in this case, the Jacobian matrix of the system evaluated at $E_{b}^{T}=(x_{b}^{T},y_{b}^{T},0)$ is given by:
\begin{eqnarray*}
    \mathbf{J}_{E_{b}^{T}} = \left[
    \begin{array}{ccc}
    1-2x_{b}^{T}-y_{b}^{T}f_{1}'(x_{b}^{T}) & -f_{1}(x_{b}^{T}) & 0 \\
    y_{b}^{T}f_{1}'(x_{b}^{T}) & 0 & -f_{2}(y_{b}^{T}) \\
    0 & 0 & 0 \\
    \end{array}\right].
\end{eqnarray*}
This Jacobian matrix has rank $2$ because $f_{1}(x_{b}^{T})$ and $f_{2}(y_{b}^{T})$ are both non-zero. This ensures that the zero is a simple eigenvalue of $\mathbf{J}_{E_{b}^{T}}$, which implies that $E_{b}^{T}$ is a non-hyperbolic equilibrium point. The eigenvectors corresponding to simple zero eigenvalues of $\mathbf{J}_{E_{b}^{T}}$ and $[\mathbf{J}_{E_{b}^{T}}]^T$ are given by: $$\mathbf{V}= \left[\begin{array}{c}
1 \\
\frac{1-2x_{b}^T-y_{b}^Tf_1'(x_{b}^T)}{f_1'(x_{b}^T)} \\
\frac{y_{b}^Tf_1'(x_{b}^T)}{f_2(y_{b}^T))}\\
\end{array}\right] ~\mbox{and}~ \mathbf{W} = \left[\begin{array}{c}
0 \\
0 \\
1\\
\end{array}\right],$$ respectively. As we have seen that the interior equilibrium point coincides with boundary equilibrium $E_b$ at $d_2=d_{2}^T$, and at the same time, $1-2x_b^T-y_b^Tf_1'(x_b^T)=1-2x_{\star}-y_{\star}f_{1}'(x_{\star})\neq0$. Using these, we obtain the following transversality conditions:
$$\mathbf{W}^T\mathbf{F}_{d_2}((x_{b}^{T},y_{b}^{T},0);d_{2}^{T}) = 0,$$
$$\mathbf{W}^T[D\mathbf{F}((x_{b}^{T},y_{b}^{T},0);d_{2}^{T})(\mathbf{V})]\,=\frac{-y_{b}^Tf_1'(x_{b}^T)}{f_2(y_{b}^T)}=\frac{-y_{b}^Tf_1'(x_{b}^T)}{d_2}\,\neq\, 0,$$ 
$$\mbox{and}~\mathbf{W}^T[D^2\mathbf{F}((x_{b}^{T},y_{b}^{T},0);d_{2}^{T})(\mathbf{V},\mathbf{V})]\,=2f_2'(y_{b}^T)\left(\frac{1-2x_{b}^T-y_{b}^Tf_1'(x_{b}^T)}{f_{1}'(x_{b}^T)}\right)\left(\frac{y_{b}^Tf_1'(x_{b}^T)}{f_2(y_{b}^T)}\right) \neq 0.$$
This implies that the system undergoes a non-degenerate transcritical bifurcation at $d_2=d_{2}^{T}$.

\subsection{Hopf bifurcation}

 Sometimes, a system shows periodic behaviour, which can be predicted through the Hopf bifurcation. We can study such Hopf-bifurcation theoretically by knowing the characteristic equation of the Jacobian matrix around the co-existing equilibrium point with the help of Liu's criterion. In this work, we have considered $d_{2}$ as the bifurcation parameter, and for this, we can write the characteristics equation of the Jacobian matrix $\mathbf{J}_{E_{*}}$ as: 
$$\lambda^3+P_2(d_{2})\lambda^2+ P_1(d_{2})\lambda+P_0(d_{2}) = 0,$$
where $P_{0}(d_{2})$, $P_{1}(d_{2})$, and $P_{2}(d_{2})$ are given at the end of Sect. \ref{EqPtSt} which have been obtained in a polynomial of $d_2$ by some manipulation. According to Liu's criterion, at the Hopf bifurcation threshold $d_{2} = d_{2}^{H}$, the functions $P_{0}(d_{2})$, $P_{1}(d_{2})$, and $P_{2}(d_{2})$ satisfy the conditions: $$P_0(d_{2})>0, P_2(d_{2})>0, \Delta (d_{2}) = 0,~\mbox{and} ~\frac{d\Delta (d_{2})}{dd_{2}} \neq 0,$$ where $\Delta (d_{2}) = P_{1}(d_{2})P_{2}(d_{2}) - P_{0}(d_{2})$.

There is no possibility of having a co-existing equilibrium point for the considered system for $d_{2}> d_{2}^{S}$. On the other hand, it can exist two co-existing equilibrium points $E_{\star}$ and $E_{\bullet}$ for $d_{2} < d_{2}^{S}$. Hence, if the system exists a Hopf bifurcation at $d_{2}^{H}$, then it has to be less than $d_{2}^{S}$. For this case, at $d_{2} = d_{2}^{H}$, the inequality $1-2x_{\star}-y_{\star}f_1'(x_{\star}) > 0$ satisfies for the equilibrium point $E_{\star}$, which violates the condition $P_{0}(d_{2}^{H})>0$. Therefore, the system does not possess a Hopf bifurcation at any values $d_{2}< d_{2}^{S}$ for the equilibrium $E_{\star}$ where the condition $y_{\star}\tilde{f}_{1}'(x_{\star}) < -1$ holds. However, there could be a possibility of Hopf bifurcation for the other equilibrium point $E_{\bullet}$. In this case, the inequality $P_{0}(d_{2}^{H}) > 0$ satisfies at $E_{\bullet}$ as $1-2x_{\bullet}-y_{\bullet}f_1'(x_{\bullet}) < 0$.

According to our assumption (b) on $\Tilde{f_2}$, we can say $\frac{d}{dy}(\Tilde{f_2}(y))|_{(x_{*},y_{*},z_{*})}=\frac{f_2'(y_{*})y_{*}-f_2(y_{*})}{y_{*}^2}<0$ which further implies $\frac{f_2(y_{*})}{y_{*}}-f_2'(y_{*})>0$. Now using the result $f(x_{*})-d_1-z_{*}f_2'(y_{*})=z_{*}(\frac{f_2(y_{*})}{y_{*}}-f_2'(y_{*}))\neq 0$, we can conclude that $\frac{d\Delta(d_2)}{dd_2}=[(z_*f_2'(y_*))\{f_1(x_*)-d_1-z_*f_2'(y_*)\}]\neq0$ for the entire range of values of $d_2$ under consideration.

Since $P_1(d_2)$ and $P_2(d_2)$ both are differentiable functions of $d_2$, by differentiating $P_1(d_2)$, we get $\frac{d}{dd_2}P_1(d_2)=z_{*}f_2'(y_{*})>0$ implies $P_1(d_2)$ is an increasing function as we increase the value of $d_2$. Again By differentiating $P_1(d_2)$, we find $\frac{d}{dd_2}(P_2(d_2))=\frac{-z_{*}}{y_{*}}<0$ which implies $P_2(d_2)$ is a decreasing function as we increase the value of $d_2$. 
 
Now if we can find some value of $d_2$ such that the conditions $P_2(d_2)>P_1(d_2)$ and $P_1(d_2)P_2(d_2)=P_0(d_2)$ holds together, then one can see the appearance of a Hopf-bifurcation at that parametric value $d_2$. Again if we can find some value of $d_2$ such that the conditions $P_2(d_2)<P_1(d_2)$ and $P_1(d_2)P_2(d_2)=P_0(d_2)$ holds together, then One can see the appearance of another Hopf-bifurcation at that parametric value $d_2$. Therefore, we can have at most two Hopf bifurcation thresholds.

\color{black}

\section{Numerical Results}{\label{NumRes}}

In this section, at first, we demonstrate the bifurcation scenario for the coexistence equilibrium point of the model \eqref{genmodel} with Holling type II functional responses by varying $d_2$ as the bifurcation parameter. The model under consideration is the classical Hastings-Powell model, our choice of parameter values are close to the values used in 
\cite{klebanoff1994A}. Fixed parameter values for the model \eqref{genmodel} with Holling type II functional responses are mentioned in the following table.


\begin{table}[h]
\begin{center}
\begin{tabular}  {|c|c|} \hline
Parameter  &  Values \\ 
\hline 
$a_1$ & 4.98 \\ 
\hline 
$b_1$ & 6.2 \\
\hline 
$a_2$ & 0.46 \\ 
\hline 
$b_2$ & 2 \\ 
\hline 
$d_1$ & 0.4 \\ 
 \hline
\end{tabular}
\caption{Parameter values for the model with Holling type II functional responses.}
\label{Table:1}
\end{center}
\end{table}

\begin{figure}[h!]
\centering
\mbox{\subfigure[]{\includegraphics[width=6cm]{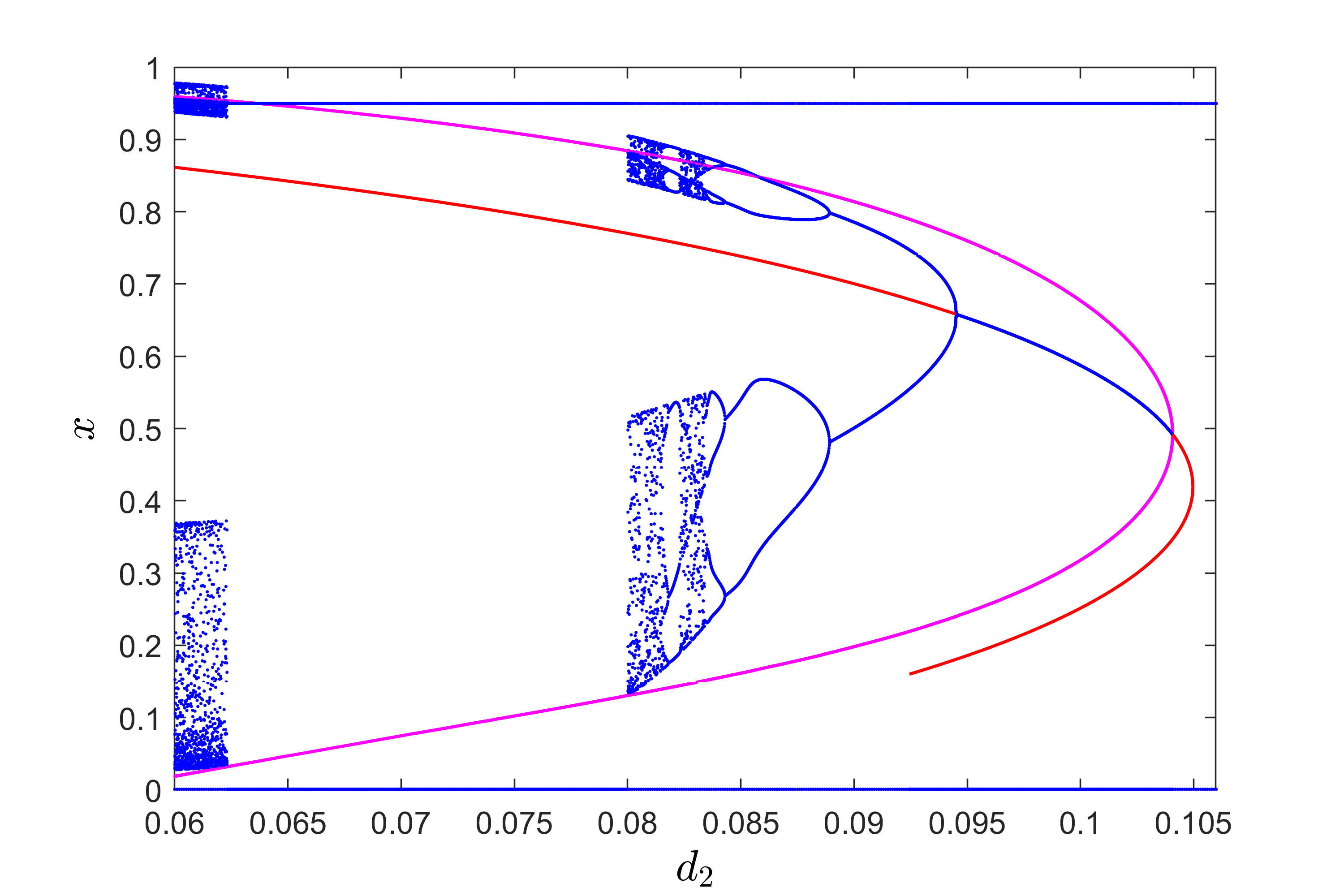}}
\subfigure[]{\includegraphics[width=6cm]{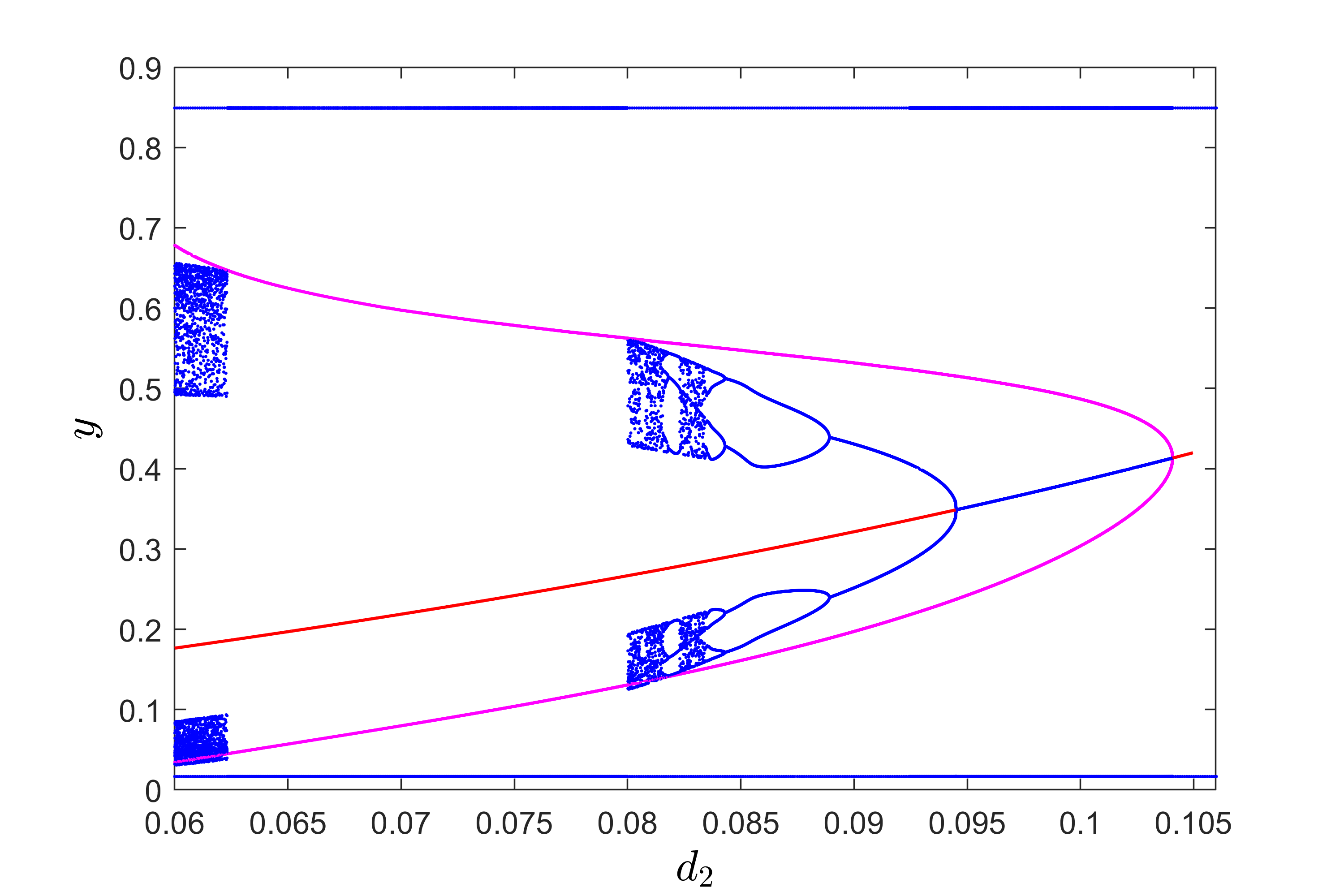}}
\subfigure[]{\includegraphics[width=6cm]{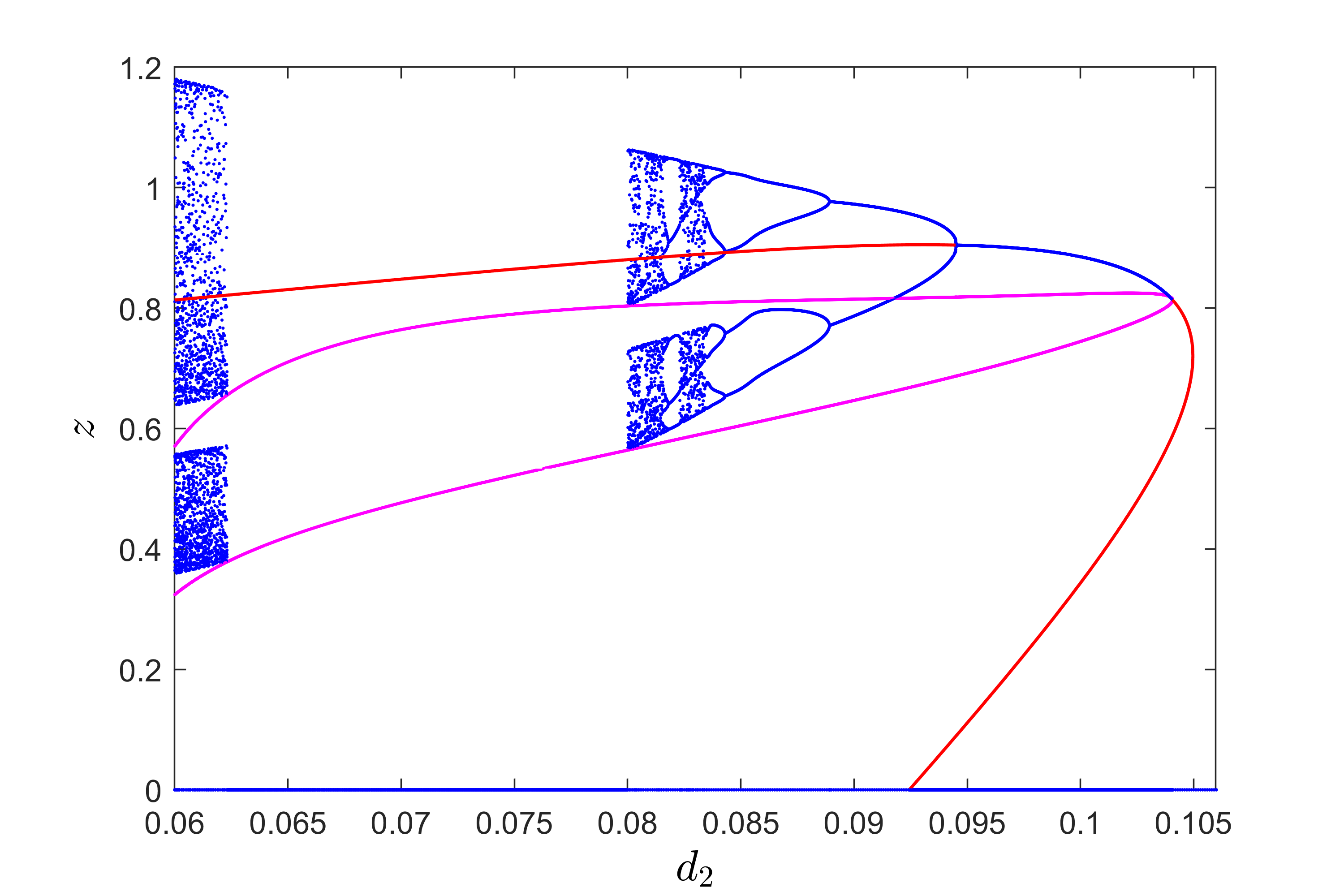}}}
\caption{Bifurcation diagram for Holling type II functional responses.}
\label{fig:fig2H}
\end{figure}

Here we will demonstrate the appearance of a chaotic attractor through a period-doubling route and its disappearance due to the collision with an unstable limit cycle arising from the coexistence equilibrium point through a subcritical Hopf bifurcation by varying our bifurcation parameter $d_2$ between the range $[0.06,0.105]$. Note that our choice of $d_1$ satisfies the condition $d_1<f_1(1)$ and hence the axial equilibrium point $E_1$ is a saddle point throughout this mentioned range of $d_2$ with $xz$ plane serving as a two-dimensional stable manifold. Based on our analytical results in the previous sections, we can calculate various bifurcation thresholds. We found a saddle-node bifurcation threshold at $d_{2}^S\equiv 0.1049651383$, two Hopf-bifurcation thresholds at $d_{2}^{H1}\equiv 0.10406993$ and $d_2^{H2}\equiv 0.09453397$ and a transcritical bifurcation threshold at $d_{2}^T\equiv 0.09244019$. Now we will discuss the dynamical behaviour of the coexistence equilibrium point(s) due to the variation of $d_2$ within the specified range. If we choose the value of $d_2$ greater than $d_{2}^S$, the system has no co-existing equilibrium point. Two co-existing equilibrium points $E_\star$ and $E_\bullet$ bifurcate through saddle-node bifurcation as the parameter $d_2$ decreases through $d_2^S$. Considering our results of the previous section, we can see that the components of these two bifurcated equilibrium points satisfy the conditions $y_{\star}\Tilde{f}_1'(x_{\star})<-1$ and $y_{\bullet}\Tilde{f}_1'(x_{\bullet})>-1$, where $y_\star=y_\bullet$. Now according to our assumptions on the properties of $\Tilde{f_1}$, we can say $y_*\Tilde{f_1'}$ is an increasing function. As a consequence, we get the inequality $x_{\star}<x_{\bullet}$ which leads to the inequality $z_{\star}<z_{\bullet}$. Therefore, once we plot the $x$ and $z$-coordinates of the interior equilibrium points in the bifurcation diagram, the lower branches will correspond to the equilibrium point $E_\star$, which has been marked in red in Figure 1(a,c). Since the $y$ components of both the equilibrium points are equal, the plot of $y$-coordinate for both the equilibrium points will coincide (see Fig.~3(b)). Due to the condition $1-2x_{\star}-y_{\star}f_1'(x_{\star}) > 0$, we can say that $E_\star$ remains unstable whenever it exists and disappears through a transcritical bifurcation when $d_2$ decreases through $d_2^T$. The other equilibrium point $E_\bullet$ is a saddle point for $d_{2}^{H1}<d_2<d_2^S$ and becomes stable as $d_2$ crosses the Hopf bifurcation threshold $d_{2}^{H1}$. An unstable limit cycle emerges at $d_{2}=d_2^{H_{sub}}$ and continues to exist for $d_{2}<d_2^{H_{sub}}$. In summary, a subcritical Hopf bifurcation occurs at $d_2=d_{2}^{H1}$. $E_\bullet$ remains stable up within the parameter range $d_2^{H2}<d_2<d_{2}^{H1}$ and losses stability through another Hopf bifurcation at $d_2=d_2^{H2}$. A stable limit cycle bifurcates from the stable equilibrium point at this threshold. This observation confirms that the Hopf bifurcation at $d_2=d_2^{H2}$ is a supercritical Hopf bifurcation. This stable limit cycle exists for $d_2<d_2^{H2}$ until it undergoes a period-doubling bifurcation. In Fig.~\ref{fig:fig2H}, stable branches of equilibrium points and limit cycles are marked in blue, components of unstable equilibrium are marked in red, and the unstable limit cycle bifurcating from the subcritical Hopf bifurcation is marked in magenta. 

The stability of the boundary equilibrium point $E_b$ depends on the parameter value $d_2$. For $d_2>d_2^T$, the boundary equilibrium point is a saddle point with two complex conjugate eigenvalues with a positive real part and one negative eigenvalue. The negative eigenvalue corresponds to a stable manifold transversal to the $xy$ plane. At the threshold value $d_2^T$, the co-existing equilibrium point $E_\star$ coincides with the boundary equilibrium point and disappears through a transcritical bifurcation. The disappearance of the co-existing equilibrium point $E_{\star}$ can be understood from the plot for the $z$ component (see Fig.~\ref{fig:fig2H}(c)); for the sake of brevity, the $x$ and $y$ components of boundary equilibrium are not shown in Fig.~\ref{fig:fig2H}(a,b). There exists a stable limit cycle on the $xy$ plane for the entire range of $d_2$ under consideration, which has been represented as two nearly horizontal lines in Fig.~\ref{fig:fig2H}(a,b). As the stable limit cycle lies on the $ xy$-plane, the plot of the $z$ has a blue line $z=0$ in Fig.~\ref{fig:fig2H}(c). Furthermore, if $d_2^S<d_2$, where $d_2^S$ is the threshold value beyond which the system has no co-existing equilibrium point, it is reasonable to conclude that the top predator goes to extinction when their death rate is significantly high, precisely for $d_2>d_2^S$. The absence of a co-existing equilibrium point suggests that the population dynamics of the top predator are not sustainable, leading to its decline and eventual extinction. The unstable limit cycle mentioned earlier acts as a separatrix between the basins of attraction of the stable boundary limit cycle (limit cycle in $xy$ plane) and the stable co-existing equilibrium point ($E_{\bullet}$), the stable limit cycle bifurcates through supercritical Hopf bifurcation as well as a chaotic attractor for the range $0.06\leq d_2\leq 0.104$ (approximately).  

The stable co-existing limit cycle undergoes period-doubling bifurcations, which involve a doubling of the period of the limit cycle as the parameter $d_2$ varies. This period-doubling cascade begins when $d_2$ falls below the value $0.089$ (approximately). This doubling occurs repeatedly, leading to increasingly complex dynamics and eventually leads to chaotic behaviour. The chaotic oscillation disappears when the chaotic attractor hits the unstable limit cycle for a value close to $d_2 = 0.08$. This collision alters the dynamics of the system drastically, leading to the cessation of chaotic behaviour. Interestingly, the chaotic behaviour reappears through another collision with the same unstable limit cycle at $d_2 \equiv 0.062$. Ecologically, this cyclical behaviour could correspond to recurring periods of instability and variability in the ecological community, with potential impacts on species coexistence, trophic interactions, and ecosystem functioning.

\subsection{Structural Sensitivity}

\begin{figure}[h!]
\centering
\mbox{\subfigure[]{\includegraphics[width=5.56cm]{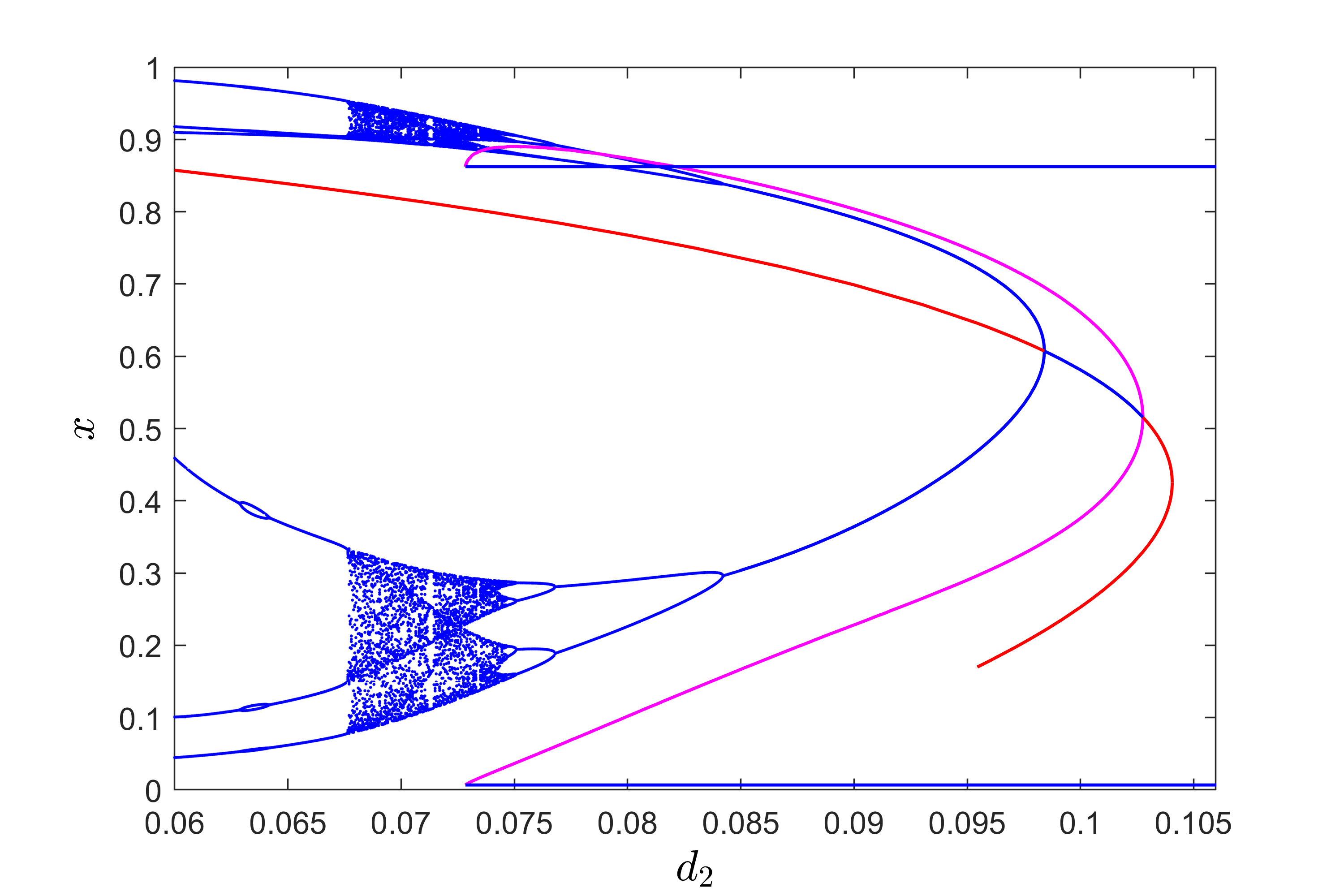}}
\subfigure[]{\includegraphics[width=5.56cm]{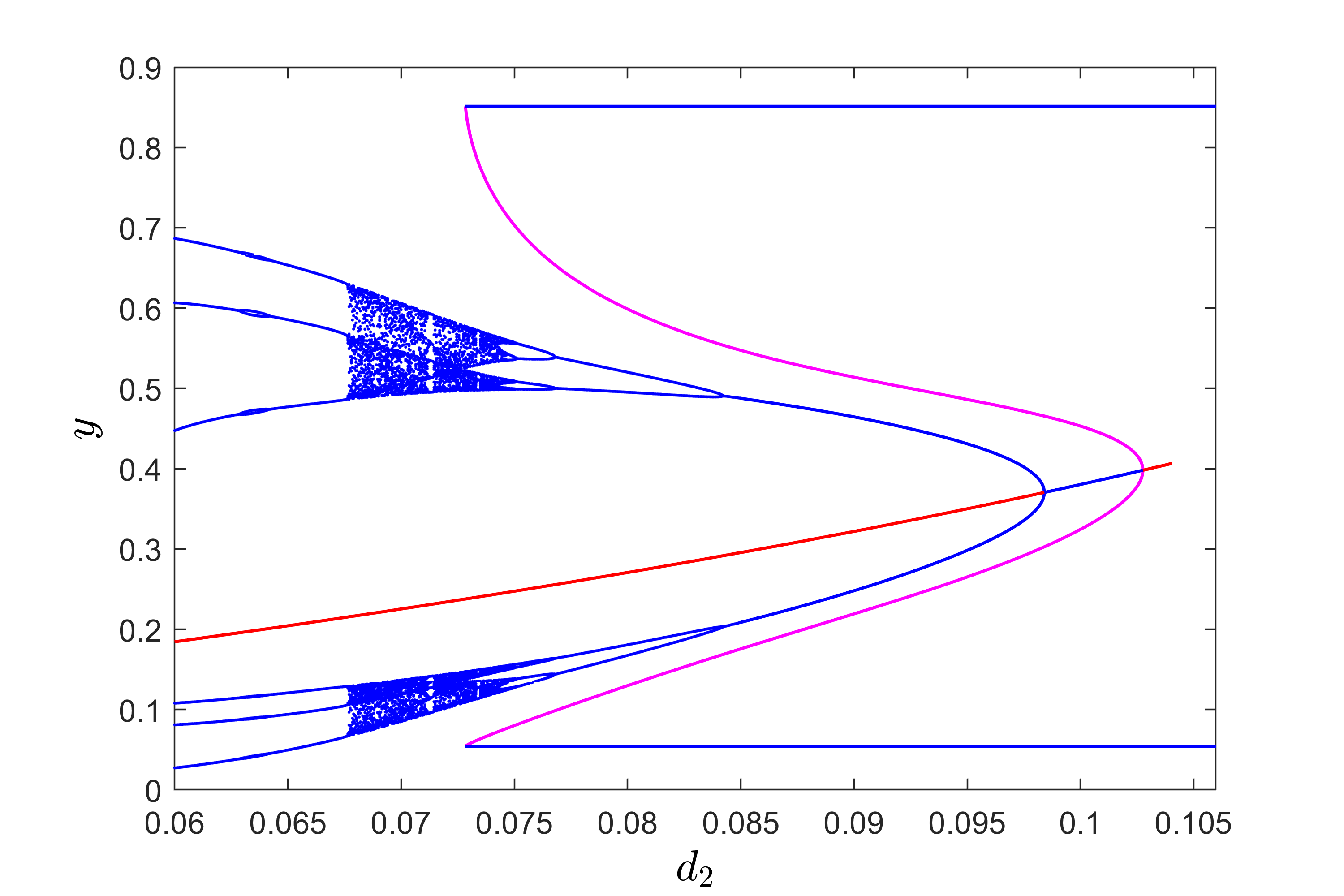}}
\subfigure[]{\includegraphics[width=5.5cm]{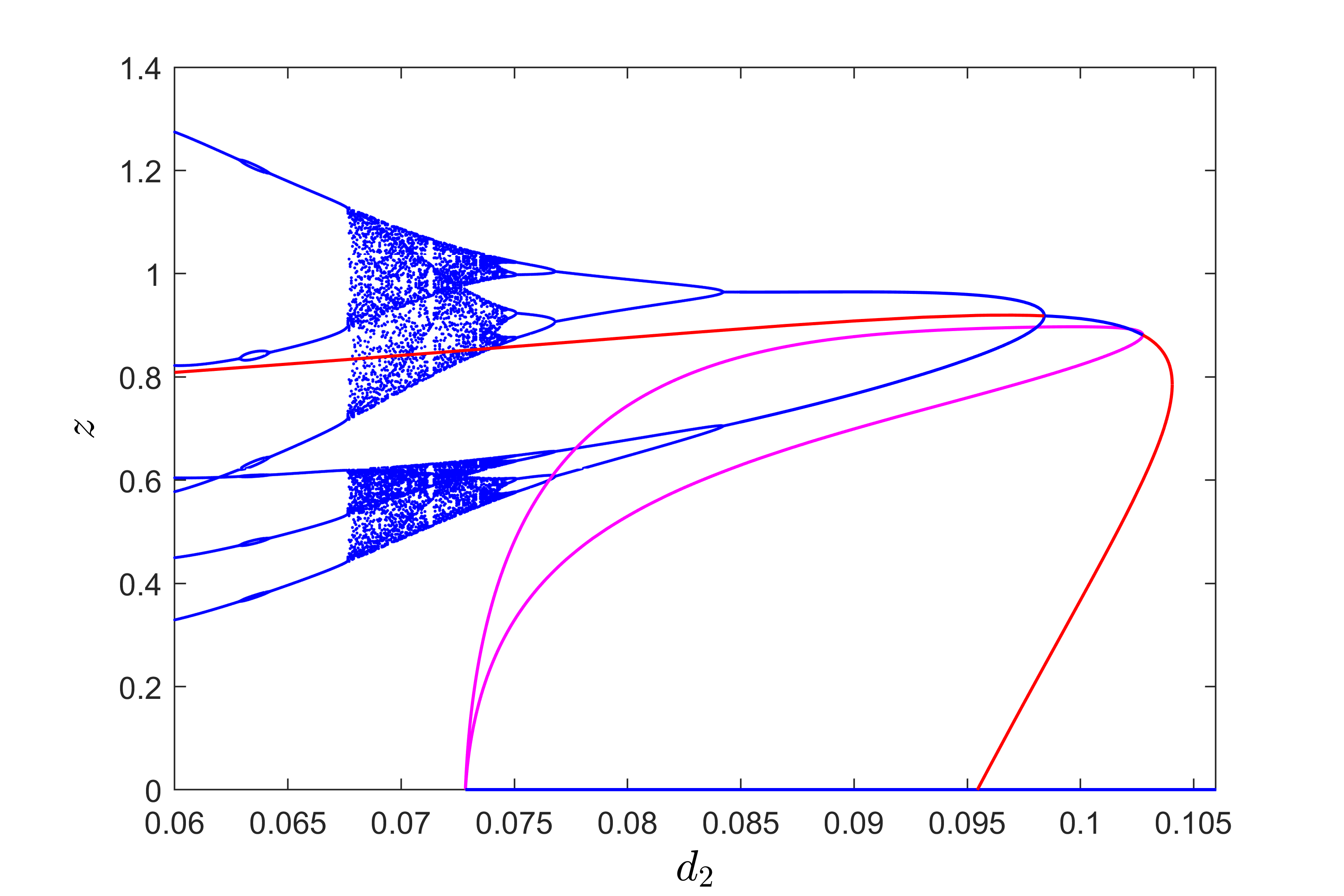}}}
\caption{Bifurcation diagram for the Ivlev functional responses.}
\label{fig:fig2I}
\end{figure}

Now we want to study the structural sensitivity of the bifurcations of the coexistence equilibrium point and the route to chaos in the model under consideration by substituting the Ivlev functional response in place of the Holling type II functional response. We will specifically focus on the parameter $d_2$ as the bifurcation parameter. To assess the structural sensitivity, we adopt the methodology outlined in the study by Fussmann et al. \cite{fussmann2005} to determine the parameter values $\bar{a}_j$ and $\bar{b}_j$ associated with the Ivlev functional responses. The fixed parameter values relevant to this transition are presented in the following table. By employing nonlinear least square regression, we determine the values of $\bar{a}_1$ and $\bar{b}_1$ while maintaining the Holling type II functional response with $a_1=4.98$ and $b_1=6.2$. Similarly, we obtain the values of $\bar{a}_2$ and $\bar{b}_2$.
\begin{table}[h]
\begin{center}
\begin{tabular}  {|l | r|} \hline
Parameter  &  Values \\ 
\hline 
$\bar{a}_1$ & 0.67 \\ 
\hline 
$\bar{b}_1$ & 5.349 \\
\hline 
$\bar{a}_2$ & 0.1647 \\ 
\hline 
$\bar{b}_2$ & 2.457 \\ 
 \hline
\end{tabular}
\caption{Parameter values for the model with Ivlev functional responses.}
\label{Table:2}
\end{center}
\end{table}
In order to simulate the model \eqref{genmodel} with Ivlev functional responses, we maintain the same value of $d_1$ as mentioned in Table-1. Other parameter values are as in Table-2 and $d_2$ as bifurcation parameters varying between the values $0.06$ and $0.105$. let us denote Ivlev functional responses as $f_{j_{Ivlev}}$ for $j=1,2$. Here also, the choice of $d_1$ satisfies the inequality $d_1<f_{1_{Ivlev}}(1)$, which ensures that the axial equilibrium point $E_{1_{Ivlev}}$ is a saddle point possessing a two-dimensional stable manifold within the $xz$ plane. The model incorporating Ivlev functional responses exhibits these bifurcation thresholds: a saddle-node bifurcation threshold at $d_{2_{Ivlev}}^S\equiv 0.10405163$, two Hopf bifurcation thresholds at $d_{2_{Ivlev}}^{H1}\equiv 0.10275556$ and $d_{2_{Ivlev}}^{H2}\equiv 0.09840295$, and a transcritical bifurcation threshold at $d_{2_{Ivlev}}^T\equiv 0.09544625$. When $d_2$ exceeds the saddle-node bifurcation threshold $d_{2_{Ivlev}}^S$, the model does not possess any co-existing equilibrium points. However, as $d_2$ decreases below $d_{2_{Ivlev}}^S$, two co-existing equilibrium points, denoted as $E_{\star_{Ivlev}}$ and $E_{\bullet_{Ivlev}}$, emerge. From previous results we can say at $E_{\star_{Ivlev}}$ and $E_{\bullet_{Ivlev}}$, we have the conditions $y_{\star}{\Tilde{f_1}}'_{Ivlev}(x_{\star})<-1$ and $y_{\bullet}{\Tilde{f_1}}'_{Ivlev}(x_{\bullet})>-1$, where $y_\star=y_\bullet$. Considering the condition $y_{\star}{\Tilde{f_1}}'_{Ivlev}(x_{\star})<-1$ and assumptions on ${\Tilde{f_1}}_{Ivlev}$, it can be deduced that $E_{\star_{Ivlev}}$ is unstable whenever it exists. It disappears through a transcritical bifurcation when $d_2$ reaches $d_{2_{Ivlev}}^T$. In Fig.~\ref{fig:fig2I}(a,c), the lower branches marked in red correspond to the equilibrium point $E_{\star_{Ivlev}}$. The other co-existing equilibrium point $E_{\bullet_{Ivlev}}$ remains a saddle point within the range $d_{2_{Ivlev}}^{H1}<d_2<d_{2_{Ivlev}}^S$ and becomes stable as $d_2$ decreases below $d_{2_{Ivlev}}^{H1}$. At this threshold, an unstable limit cycle bifurcates and persists for $d_2<d_{2_{Ivlev}}^{H1}$, which suggests that we have a subcritical Hopf bifurcation at $d_2=d_{2_{Ivlev}}^{H1}$.  $E_{\bullet_{Ivlev}}$ remains stable until $d_2$ crosses the threshold $d_{2_{Ivlev}}^{H2}$, after which it loses stability as $d_2$ further decreases. At this threshold, a stable limit cycle emerges, which exists until it undergoes a period-doubling bifurcation. This fact says that we have a supercritical Hopf bifurcation at $d_2=d_{2_{Ivlev}}^{H2}$. The overall behaviour can be observed in the upper branch of Fig.~\ref{fig:fig2I}(a,c).

The model exhibits a boundary equilibrium point and a stable boundary limit cycle in the $xy$ plane for the considered values of $d_2$. The co-existing equilibrium point $E_{\star_{Ivlev}}$ collides with the boundary equilibrium point at $d_2=d_{2_{Ivlev}}^T$, causing a transcritical bifurcation. The stable co-existing limit cycle undergoes a period-doubling bifurcation, leading to chaotic behaviour, as $d_2$ decreases through the value $d_2=0.084$. The chaotic attractor persists within the range $[0.0677,0.0745]$ and did not hit the unstable limit cycle generated from the subcritical Hopf bifurcation. Whereas, in this scenario,  the unstable limit cycle disappears at $d_2\equiv 0.073$ through a global bifurcation when it collides with the stable boundary limit cycle (limit cycle in $xy$ plane). The chaotic attractor continues to exist over an extended range of parameter values and eventually ceases to exist through a crisis, leading to the emergence of a three-periodic coexistence scenario. The disappearance of the unstable limit cycle suggests a significant shift in the population dynamics, potentially leading to a different ecological regime. This indicates the robustness and resilience of the chaotic dynamics in the ecological system. However, eventually, the chaotic attractor ceases to exist through a crisis, which represents a sudden and drastic change in the population dynamics. This crisis event leads to the emergence of a three-periodic coexistence scenario, where the interacting species exhibit cyclical patterns of population fluctuations with three distinct periods.

In summary, the observed dynamics highlight the complexity and sensitivity of population interactions in the ecological system. The coexistence of chaotic attractors, the disappearance of unstable limit cycles, and the emergence of periodic coexistence scenarios all reflect the intricate interplay between species interactions, ecological conditions, and the parameter values governing the system.

Here we want to mention that the limit cycle lying on the $xy$ plane is plotted in Fig.~\ref{fig:fig2H} for the entire range of parameters for Holling type II functional responses. In this case, we find a non-empty basin of attraction for the stable limit cycle on $xy$-plane in the interior of the first octant.
However, there is no initial condition in the interior of the first octant from which we can reach the stable limit cycle on the $xy$-plane for the model with Ivlev functional responses when $d_2<0.073$. As a matter of fact, the branches of the limit cycle lying on the $xy$-plane are not shown in Fig.~\ref{fig:fig2I} for $d_2<0.073$. 

\vspace{1cm}

\section{Conclusions}

\begin{figure}[h!]
\centering
\mbox{\subfigure[]{\includegraphics[width=8cm]{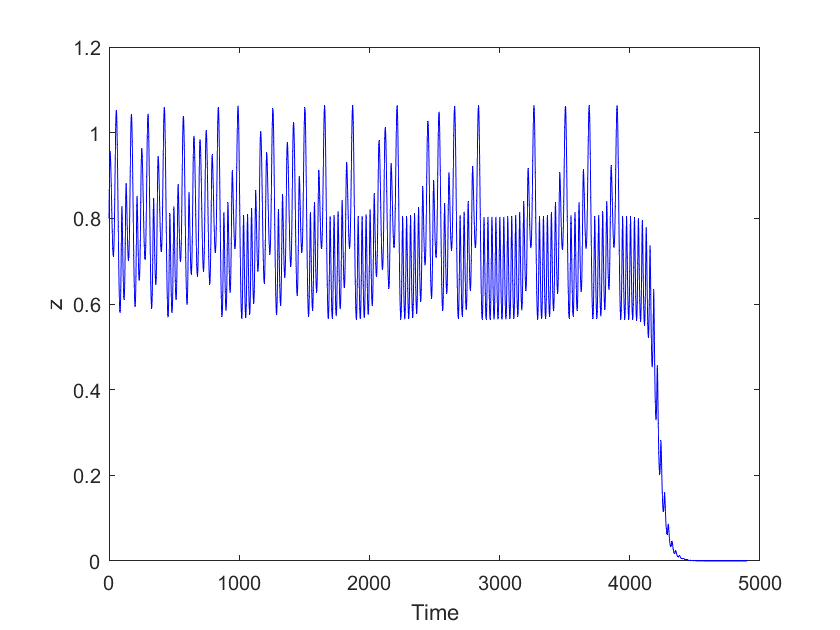}}
\subfigure[]{\includegraphics[width=8cm]{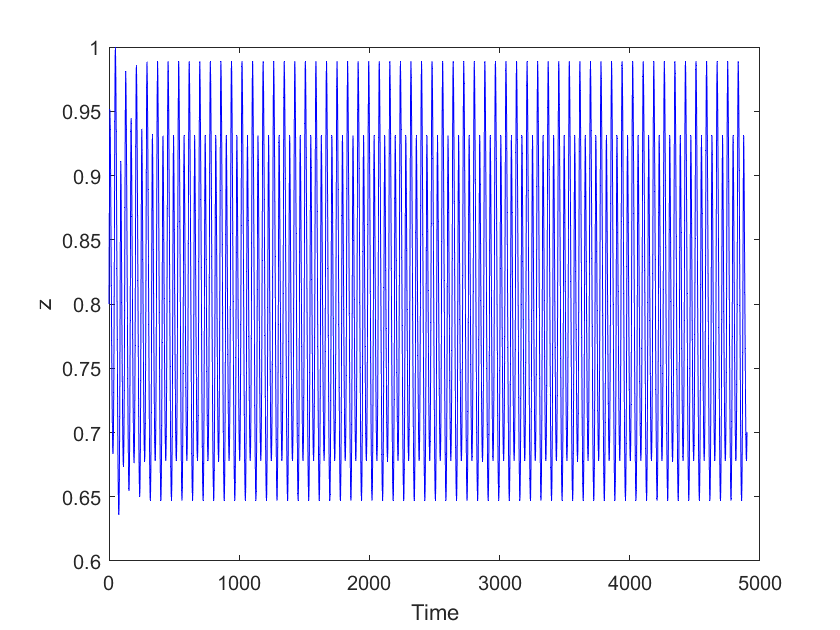}}}
\caption{ The time evolution of top predator density $z$ as (a) for the Holling type II functional responses (with parameter values $a_1=4.98,\,b_1=6.2,\,a_2=0.46,\,b_2=2,\,d_1=0.4,\,d_2=0.08$), (b) for the Ivlev functional responses (with parameter values $\bar{a}_1=0.67,\,\bar{b}_1=5.349,\,\bar{a}_2=0.1647,\,\bar{b}_2=2.457,\,d_1=0.4,\,d_2=0.08$).}
\label{fig:fig5}
\end{figure}

Hastings-Powell model is well known as the first prototype model for prey-predator-top-predator type interactions which can explain the complex oscillatory coexistence of three interacting populations. In the classical Rosenzweig-MacArthur model, we find stable oscillatory coexistence with high amplitude oscillation for both the prey and predator species beyond the supercritical Hopf bifurcation threshold. However, a large stable limit cycle situated close to the coordinate axes indicates the possibility of extinction of one or more species under environmental and/or demographic variability \cite{may1976B}. Nowadays, researchers are interested in transient dynamics and the emergence of extinction scenarios through global bifurcations apart from the structural sensitivity of the models under consideration. The Hastings-Powell model is capable of producing long transients before the extinction of top predators, as shown in Fig.~\ref{fig:fig5}. This kind of dynamics was first reported by Yodzis in \cite{mccann1994B}. The existence of long transients and then the extinction of the top predator is presented for the parameter value close to $d_2=0.08$ where the chaotic attractor hits the unstable limit cycle and disappears. This phenomenon is structurally sensible for the model under consideration.

To justify the structural sensitivity of the long transients leading to the extinction of top predator, we consider the model \eqref{genmodel} with Ivlev functional responses, $d_2=0.08$ and other parameter values as mentioned in Table 2 in the previous section, with initial condition ($0.45,0.5,0.8$). In this situation, we find the coexistence scenario for all three constituent species as the chaotic attractor remains quite a bit away from the unstable limit cycle. The unstable limit cycle emerged through the sub-critical Hopf bifurcation and exists for the range $d_2\in[0.073,0.1027]$ approximately. The unstable limit cycle disappears through a global bifurcation when it collides with the stable boundary limit cycle (limit cycle on the $xy$-plane). Comparing the bifurcation diagrams in 
Figs.~\ref{fig:fig2H} and \ref{fig:fig2I}, we conclude that the survival and extinction scenarios are solely influenced by the parametrization of the functional responses. Alternative parametrization of functional responses, without changing the basic properties, indicates how small changes in the population of one species affect the growth of other species, which are directly or indirectly related through the functional responses. This finding indicates that the expected extinction of one or more species may be avoided due to some adaptation mechanism which is reflected through the change in the parametrization of the functional responses.

In addition, considering alternative parametrization for the functional responses does not alter the overall bifurcation scenario. However, the size of the chaotic attractor and position of the chaotic attractor with respect to the position of the unstable limit cycle is sensitive with respect to the parametrization of functional responses. The model with Holling type II functional responses predicts the extinction of the top predator when $d_2$ lies in the range [$0.062,0.08$]. There is no possibility for the extinction of top predators if we replace the Holling type II functional responses with Ivlev functional responses. 

\bibliography{References}

\end{document}